\documentclass[reqno,11pt]{amsart}
\usepackage{amscd,amssymb,verbatim}
\numberwithin{equation}{section} \theoremstyle{plain}

\newcommand\alp{\alpha}         
\newcommand\bet{\beta}
\newcommand\gam{\gamma}         \newcommand\Gam{\Gamma}
         
\newcommand\eps{\varepsilon}
\newcommand\zet{\zeta}
\newcommand\tet{\theta}

\newcommand\lam{\lambda}                \newcommand\Lam{\Lambda}
\newcommand\sig{\sigma}         \newcommand\Sig{\Sigma}

\newcommand\ome{\omega}         \newcommand\Ome{\Omega}

\newcommand\calB{{\mathcal{B}}}
\newcommand\calC{{\mathcal{C}}}
\newcommand\calD{{\mathcal{D}}}

\newcommand\calI{{\mathcal{I}}}

\newcommand\calN{{\mathcal{N}}}
\newcommand\calO{{\mathcal{O}}}

\newcommand\calR{{\mathcal{R}}}

\newcommand\bfw{{\mathbf w}}

\newcommand\RR{\mathbb{R}}

\newcommand\ZZ{\mathbb{Z}}

\newcommand\CC{\mathbb{C}}

 \newcommand\gro{{\mathfrak{o}}}

\newcommand\nek{,\ldots,}
\newcommand\sdp{\times \hskip -0.3em {\raise 0.3ex
\hbox{$\scriptscriptstyle |$}}} 

\newcommand\Cone{\operatorname{Cone}}

\newcommand\Det{\operatorname{Det}}

\newcommand\End{\operatorname{End\,}}

\newcommand\IM{\operatorname{Im}}

\newcommand\MOD{\operatorname{mod}}

\newcommand\Ker{\operatorname{Ker}}

\newcommand\rank{\operatorname{rank}}

\newcommand\RE{\operatorname{Re}}

\newcommand{\sign}{\operatorname{sign}}

\newcommand\Tr{\operatorname{Tr}}

\newcommand\hatL{{\widehat{L}}}

\newcommand\tilE{{\widetilde{E}}}\newcommand\tile{{\widetilde{e}}}

\newcommand\tilM{{\widetilde{M}}}

\newcommand\tilw{{\tilde{w}}}

\newcommand\tilx{{\tilde{x}}}

\newcommand\tilGam{{\widetilde{\Gamma}}}

\renewcommand{\>}{\rangle}
\newcommand{\<}{\langle}

\theoremstyle{plain}
\newtheorem{Thm}[subsection]{Theorem}
\newtheorem{Cor}[subsection]{Corollary}
\newtheorem{Lem}[subsection]{Lemma}
\newtheorem{Prop}[subsection]{Proposition}
\newtheorem{Conjec}[subsection]{Conjecture}

\newtheorem{Def}[subsection]{Definition}

\theoremstyle{remark}

\newtheorem{Rem}[subsection]{Remark}

\def\TeXref#1{%
        \leavevmode\vadjust{\setbox0=\hbox{{\tt
                \  {\tiny \textrm #1}}}%
        \theight=\ht0
        \advance\theight by \lineskip
        \kern -\theight \vbox to
        \theight{\rightline{\rlap{\box0}}%
        \vss}%
        }}%
\newif\ifShowLabels
\ShowLabelstrue
\newdimen\theight
\def\TeXrefEq#1{%
        \leavevmode\vadjust{\setbox0=\hbox{{\tt
                \  {\tiny \textrm #1}}}%
        \theight=\ht1
        \advance\theight by \lineskip
        \kern -\theight \vbox to
        \theight{\rightline{\rlap{\box0}}%
        \vss}%
        }}%

\newcommand{\refs}[1]{Section ~\ref{S:#1}}
\newcommand{\refss}[1]{Subsection ~\ref{SS:#1}}

\newcommand{\reft}[1]{Theorem ~\ref{T:#1}}
\newcommand{\refl}[1]{Lemma ~\ref{L:#1}}
\newcommand{\refp}[1]{Proposition ~\ref{P:#1}}
\newcommand{\refc}[1]{Corollary ~\ref{C:#1}}

\newcommand{\refd}[1]{Definition ~\ref{D:#1}}
\newcommand{\refr}[1]{Remark ~\ref{R:#1}}
\newcommand{\refe}[1]{\eqref{E:#1}}

\newenvironment{thm}[1]%
        { \begin{Thm} \label{T:#1}  \ifShowLabels \TeXref{T:#1} \fi }%
        { \end{Thm} }

\renewcommand{\th}[1]{\begin{thm}{#1}  }
\renewcommand{\eth}{\end{thm} }

\newenvironment{lemma}[1]%
        { \begin{Lem} \label{L:#1}  \ifShowLabels \TeXref{L:#1} \fi }%
        { \end{Lem} }

\newcommand{\lem}[1]{\begin{lemma}{#1} }
\newcommand{\elem}{\end{lemma}}

\newenvironment{propos}[1]%
        { \begin{Prop} \label{P:#1}  \ifShowLabels \TeXref{P:#1} \fi }%
        { \end{Prop} }

\newcommand{\prop}[1]{\begin{propos}{#1} }
\newcommand{\eprop}{\end{propos}}

\newenvironment{corol}[1]%
        { \begin{Cor} \label{C:#1}  \ifShowLabels \TeXref{C:#1} \fi }%
        { \end{Cor} }
\newcommand{\cor}[1]{\begin{corol}{#1}  }
\newcommand{\ecor}{\end{corol}}

\newenvironment{conjec}[1]%
        { \begin{Conjec} \label{Conj:#1}  \ifShowLabels \TeXref{C:#1} \fi }%
        { \end{Conjec} }
\newcommand{\conj}[1]{\begin{conjec}{#1}  }
\newcommand{\econj}{\end{conjec}}

\newenvironment{defeni}[1]%
        { \begin{Def} \label{D:#1}  \ifShowLabels \TeXref{D:#1} \fi }%
        { \end{Def} }
\newcommand{\defe}[1]{\begin{defeni}{#1}  }
\newcommand{\edefe}{\end{defeni}}

\newenvironment{remark}[1]%
        { \begin{Rem} \label{R:#1}  \ifShowLabels \TeXref{R:#1} \fi }%
        { \end{Rem} }
\newcommand{\rem}[1]{\begin{remark}{#1}}
\newcommand{\erem}{\end{remark}}

\newcommand{\eq}[1]%
        { \ifShowLabels \TeXrefEq{E:#1} \fi
           \begin{equation} \label{E:#1} }
\newcommand{\eeq}{\end{equation}}
\newcommand{\meq}[1]%
        { \ifShowLabels \TeXrefEq{E:#1} \fi
           \begin{multline} \label{E:#1} }
\newcommand{\emeq}{\end{multline}}

\newcommand{\prf}{ \begin{proof} }
\newcommand{\eprf}{ \end{proof} }
\newcommand{\Label}[1]{\label{#1}  \ifShowLabels \TeXref{#1} \fi }

\ShowLabelsfalse

\newcommand{\n}{\nabla}
\newcommand{\na}{\n_\alp}\newcommand{\Ea}{E_\alp}

\renewcommand{\b}{\bullet}
\newcommand{\pa}{\text{\( \partial\)}}\newcommand{\tilpa}{\text{\( \tilde{\partial}\)}}\newcommand{\hatpa}{\text{\( \hat{\partial}\)}}

\newcommand{\even}{{\operatorname{even}}}

\newcommand{\Arg}{\operatorname{\mathbf{Arg}}}
\newcommand{\Ph}{\operatorname{\mathbf{Ph}}}
\newcommand{\Mon}{\operatorname{Mon}}

\newcommand{\rabs}{\rho^{\operatorname{abs}}}

\newcommand{\Detgrtetnp}{\Det_{{\operatorname{gr},\tet}}}

\newcommand{\p}{\pi_1(M)}
\newcommand{\Rep}{\operatorname{Rep}(\p,\CC^n)}
\newcommand{\Repo}{\operatorname{Rep_0}(\p,\CC^n)}

\newcommand{\Repho}{\operatorname{Rep}^u_0(\p,\CC^n)}

\newcommand{\ccomp}{\calC}

\newcommand{\B}{\calB}

\newcommand{\rat}{\rho_{\operatorname{an}}}
\newcommand{\rtt}{\rho_{\operatorname{\eps,\gro}}}
\newcommand{\RSn}[1]{\|#1\|_{\Det(H^\b(M,E_\alp))}^{\operatorname{RS}}}

\newcommand{\Dbundle}{\calD{}et}
\newcommand{\Mat}{\operatorname{Mat}}
\newcommand{\trivial}{{\operatorname{trivial}}}
\begin{document}

\title{Ray-Singer Type Theorem for the Refined Analytic Torsion}
\author[Maxim Braverman]{Maxim Braverman$^\dag$}
\address{Department of Mathematics\\
        Northeastern University   \\
        Boston, MA 02115 \\
        USA
         }
\email{maximbraverman@neu.edu}
\author[Thomas Kappeler]{Thomas Kappeler$^\ddag$}
\address{Institut fur Mathematik\\
         Universitat Z\"urich\\
         Winterthurerstrasse 190\\
         CH-8057 Z\"urich\\
         Switzerland
         }
\email{tk@math.unizh.ch}
\thanks{${}^\dag$Supported in part by the NSF grant DMS-0204421.\\
\indent${}^\ddag$Supported in part by the Swiss National Science foundation, the programme SPECT, and the European Community through the FP6
Marie Curie RTN ENIGMA (MRTN-CT-2004-5652)}
\begin{abstract}
We show that the refined analytic torsion is a holomorphic section of the determinant line bundle over the space of complex representations of
the fundamental group of a closed oriented odd dimensional manifold. Further, we calculate the ratio of the refined analytic torsion and the
Farber-Turaev combinatorial torsion.

As an application, we establish a formula relating the eta-invariant and the phase of the Farber-Turaev torsion, which extends a theorem of
Farber and earlier results of ours. This formula allows to study the spectral flow using methods of combinatorial topology.
\end{abstract}
\maketitle

\section{Introduction}\Label{S:introduction}

Let $M$ be a closed oriented odd dimensional manifold. Denote by $\Rep$ the space of $n$-dimensional complex representations of the fundamental group
$\p$ of $M$. For $\alp\in \Rep$ we denote by $\Ea$ the flat vector bundle over $M$ whose monodromy is equal to $\alp$. Let $\na$ be the flat
connection on $\Ea$.  In \cite{BrKappelerRATdetline}, we defined the non-zero element
\[
    \rat(\alp)\ = \ \rat(\na) \ \in \ \Det\big(\, H^\b(M,\Ea)\,\big)
\]
of the determinant line $\Det\big(H^\b(M,\Ea)\big)$ of the cohomology $H^\b(M,\Ea)$ of $M$ with coefficients in $\Ea$. This element, called the
{\em refined analytic torsion}, carries information about the Ray-Singer metric and about the $\eta$-invariant. In particular, if $\alp$ is a
unitary representation, then the Ray-Singer norm of $\rat(\alp)$ is equal to 1.

\subsubsection*{Analyticity of the refined analytic torsion}
The disjoint union of the lines $\Det\big(H^\b(M,\Ea)\big)$, $(\alp\in \Rep)$, forms a line bundle $\Dbundle\to \Rep$, called the {\em
determinant line bundle}, cf. \cite[\S9.7]{BeGeVe}. It admits a nowhere vanishing section, given by the Farber-Turaev torsion, and, hence, has a
natural structure of a trivializable holomorphic bundle.

Our first result is that $\rat(\alp)$ is a nowhere vanishing holomorphic section of the bundle $\Dbundle$. It means that the ratio of the
refined analytic and the Farber-Turaev torsions is a holomorphic function on $\Rep$. For an acyclic representation $\alp$, the determinant line
$\Det\big(H^\b(M,\Ea)\big)$ is canonically isomorphic to $\CC$ and $\rat(\alp)$ can be viewed as a non-zero complex number. We show that
$\rat(\alp)$ is a holomorphic function on the open set $\Repo\subset \Rep$ of acyclic representations. This result extends Corollary~13.11 of
\cite{BrKappelerRAT}. See also \cite{BurgheleaHaller_function} for somewhat related results.

\subsubsection*{Comparison with the Farber-Turaev torsion}
In \cite{Turaev86,Turaev90}, Turaev constructed a refined version of the combinatorial torsion associated to an acyclic representation $\alp$, which depends on
additional combinatorial data, denoted by $\eps$ and called the {\em Euler structure}, as well as on the {\em cohomological orientation} of $M$, i.e., on the
orientation $\gro$ of the determinant line of the cohomology $H^\b(M,\RR)$ of $M$. In \cite{FarberTuraev00}, Farber and Turaev extended the definition of the Turaev
torsion to non-acyclic representations. The Farber-Turaev torsion associated to a representation $\alp$, an Euler structure $\eps$, and a cohomological orientation
$\gro$ is a non-zero element $\rtt(\alp)$ of the determinant line $\Det\big(H^\b(M,\Ea)\big)$.


\reft{Ray-Singer} of this paper states, that for each connected component%
\footnote{In this paper we always consider the classical (not the Zariski) topology on the complex analytic space $\Rep$.}
$\calC$ of the space\linebreak $\Rep$, there exists a constant $\tet\in \RR$, such that%
\footnote{Note that since $\rat(\alp)$ and $\rtt(\alp)$ are non-vanishing sections of the same bundle, their ratio is a non-zero complex-valued function.}
\eq{intratio}
    \frac{\rat(\alp)}{\rtt(\alp)} \ = \ e^{i\tet}\cdot f_{\eps,\gro}(\alp),
\end{equation}
where $f_{\eps,\gro}(\alp)$ is a holomorphic function of $\alp\in \Rep$, given by an explicit local expression, cf. \refe{logTTTur}. In the case
where $\alp$ is an acyclic representation close to an acyclic unitary representation, this formula was obtained in \cite{BrKappelerRAT}.

Recently, Rung-Tzung~Huang \cite{Huang06} showed by an explicit calculation for lens spaces that the constant $\tet$ can depend on the connected
component $\ccomp$.

\subsubsection*{Sketch of the proof of formula \refe{intratio}}
Using the calculation of the Ray-Singer norm of the Farber-Turaev torsion, given in Theorem~10.2 of \cite{FarberTuraev00} and the formula for
the Ray-Singer norm of the refined analytic torsion \cite[Th.~11.3]{BrKappelerRATdetline}, we obtain (cf. \refe{normofration2}) that
\eq{intrationav}
    \left|\,\frac{\rat(\alp)}{\rtt(\alp)}\,\right| \ = \  |f_{\eps,\gro}(\alp)|.
\end{equation}
Both, the left and the right hand side of this equality, are absolute values of holomorphic functions. If the absolute values of two holomorphic
functions are equal, then the two functions are equal up to a multiplication by a locally constant function, whose absolute value is equal to one.
Hence, \refe{intratio} follows from \refe{intrationav}.

\subsubsection*{Application: relation of the $\eta$-invariant with the phase of the Farber-Turaev torsion}
If $\alp\in \linebreak\Repo$ is an acyclic unitary representation, then the refined analytic torsion $\rat(\alp)$ is a non-zero complex number,
whose phase is equal, up to a correction term, to the $\eta$-invariant $\eta_\alp$ of the odd signature operator corresponding to the flat
connection on $\Ea$, cf. \refe{PhT0}. Hence, if $\alp_1$ and $\alp_2$ are two acyclic unitary representations which lie in the same connected
component of $\Rep$, equality  \refe{intratio} allows to compute the difference $\eta_{\alp_1}-\eta_{\alp_2}$ in terms of the phases of the
Farber-Turaev torsions $\rtt(\alp_1)$ and $\rtt(\alp_2)$. The significance of this computation is that it allows to study the spectral invariant
$\eta_\alp$ by the methods of combinatorial topology. With some additional assumptions on $\alp_1$ and $\alp_2$ a similar result was established
in \cite{Farber00AT} and \cite{BrKappelerRAT}, cf. \refr{Farber}.

\subsubsection*{Related works}
In \cite{Turaev86,Turaev90}, Turaev  constructed a refined version of the combinatorial torsion and posed the problem of constructing its analytic analogue. In
\cite[\S10.3]{FarberTuraev00}, Farber and Turaev asked this question in a more general setting and also suggested that such an analogue should involve the
$\eta$-invariant. The proposed notion of refined torsion gives an affirmative answer to this question in full generality.

Having applications in topology in mind,  quite some time ago, Burghelea asked the question if there exists a holomorphic function on the space of acyclic
representations $\Repo$ whose absolute value is equal to the (modified) Ray-Singer torsion. In \cite{BurgheleaHaller_Euler,BurgheleaHaller_function}, Burghelea and
Haller constructed such a holomorphic function. In particular, in \cite{BurgheleaHaller_function} they outlined a construction of this function involving Laplace-type
operators acting on forms.%
\footnote{Added in proof: For a more detailed presentation see \cite{BurgheleaHaller_function2}.} %
They require that the given complex vector bundle admits a non-degenerate symmetric bilinear form, which they use to define their Laplace-type operators. The function
constructed in \cite{BurgheleaHaller_function} is similar to the invariant $\xi$ defined in \S7 of our paper \cite{BrKappelerRAT}. Burghelea and Haller then express
the square of the Farber-Turaev torsion in terms of these determinants and some additional ingredients. Hence they obtain a formula for the Farber-Turaev torsion in
terms of analytic quantities up to a sign. This result should be compared with our formula \refe{intratio}, which expresses the Farber-Turaev torsion including its
sign in analytic terms. The sign is important, in particular, for the application discussed in \refs{Fraberth}. Note that the result of Burghelea and Haller is valid
on a manifold of arbitrary, not necessarily odd dimension. Their holomorphic function is different from our refined analytic torsion and is not related to the
Atiyah-Patodi-Singer $\eta$-invariant. In \cite{BrKappelerRATdetline_BH}, we obtain an explicit formula computing the Burghelea-Haller torsion in terms of the refined
analytic torsion and the $\eta$-invariant.

\section{The Refined Analytic Torsion}\label{S:RAT}

In this section we recall the definition of the refined analytic torsion from \cite{BrKappelerRATdetline}. The refined analytic torsion is
constructed in 3 steps: first, we define the notion of refined torsion of a finite dimensional complex endowed with a chirality operator, cf.
\refd{refinedtorsion}. Then we  fix a Riemannian metric $g^M$ on $M$ and consider the odd signature operator $\B = \B(\n,g^M)$ associated to a flat
vector bundle $(E,\n)$, cf. \refd{oddsign}. Using the {\em graded determinant} of $\B$ and the definition of the refined torsion of a finite
dimensional complex with a chirality operator we construct an element $\rho= \rho(\n,g^M)$ in the determinant line of the cohomology, cf. \refe{rho}.
The element $\rho$ is almost the refined analytic torsion. However, it might depend on the Riemannian metric $g^M$ (though it does not if
$\dim{}M\equiv1\ (\MOD 4)$ or if $\rank(E)$ is divisible by 4). Finally we ``correct" $\rho$ by multiplying it by an explicit factor, the metric
anomaly of $\rho$, to obtain a diffeomorphism invariant $\rat(\n)$ of the triple $(M,E,\n)$, cf. \refd{rat}.

\subsection{The determinant line of a complex}\label{S:detline}
Given a complex vector space $V$ of dimension $\dim{}V= n$, the {\em determinant line} of $V$ is the line $\Det(V):= \Lam^nV$, where $\Lam^nV$
denotes the $n$-th exterior power of $V$. By definition, we set $\Det(0):= \CC$. Further, we denote by $\Det(V)^{-1}$ the dual line of $\Det(V)$.

Let
\[
    \begin{CD}
       (C^\b,\pa):\quad  0 \ \to C^0 @>{\pa}>> C^1 @>{\pa}>>\cdots @>{\pa}>> C^d \ \to \ 0
    \end{CD}
\]
be a  complex of finite dimensional complex vector spaces. We call the integer $d$ the {\em length} of the complex $(C^\b,\pa)$ and we denote by
$H^\b(\pa)=\bigoplus_{i=0}^d H^{i}(\pa)$ the cohomology of $(C^{\b},\pa)$. Set
\eq{DetCH}
    \Det(C^\b) \ := \ \bigotimes_{j=0}^d\,\Det(C^j)^{(-1)^j}, \qquad \Det(H^\b(\pa))\ := \ \bigotimes_{j=0}^d\,\Det(H^j(\pa))^{(-1)^j}.
\end{equation}
The lines $\Det(C^\b)$ and $\Det(H^\b(\pa))$ are referred to as the {\em determinant line of the complex} $C^\b$ and the {\em determinant line of its
cohomology}, respectively. There is a canonical isomorphism
\eq{isomorphism}
     \phi_{C^\b} \ = \ \phi_{(C^\b,\pa)}:\, \Det(C^\b) \ \longrightarrow \ \Det(H^\b(\pa)),
\end{equation}
cf., for example, \S2.4 of \cite{BrKappelerRATdetline}.

\subsection{The refined torsion of a finite dimensional complex with a chirality operator}\label{SS:rtfd}
Let $d= 2r-1$ be an odd integer and let $(C^\b,\pa)$ be a length $d$ complex of finite dimensional complex vector spaces. A {\em chirality operator}
is an involution $\Gam:C^\b\to C^\b$ such that $\Gam(C^j)= C^{d-j}$, $j=0\nek d$. For $c_j\in \Det(C^j)$ $(j=0\nek d)$ we denote by $\Gam{}c_j\in
\Det(C^{d-j})$ the image of $c_j$ under the isomorphism $\Det(C^j)\to \Det(C^{d-j})$ induced by $\Gam$.

Fix non-zero elements $c_j\in \Det(C^j)$, $j=0\nek r-1$ and denote by $c_j^{-1}$ the unique element of $\Det(C^j)^{-1}$ such that $c^{-1}_j(c_j)=1$.
Consider the element
\eq{Gamc}
    c_{{}_\Gam} \ := \ (-1)^{\calR(C^\b)}\cdot
    c_0\otimes c_1^{-1}\otimes \cdots \otimes c_{r-1}^{(-1)^{r-1}}\otimes (\Gam c_{r-1})^{(-1)^r}\otimes (\Gam c_{r-2})^{(-1)^{r-1}}
    \otimes \cdots\otimes (\Gam c_0)^{-1}
\end{equation}
of $\Det(C^\b)$, where

\eq{R(C)}
  \calR(C^\b) \ := \ \frac12\ \sum_{j=0}^{r-1}\, \dim C^{j}\cdot\big(\, \dim C^j+(-1)^{r+j}\,\big).
\end{equation}
It follows from the definition of $c_j^{-1}$ that $c_{{}_\Gam}$ is independent of the choice of $c_j$ ($j=0\nek r-1$).
\defe{refinedtorsion}
The {\em refined torsion} of the pair $(C^\b,\Gam)$ is the element
\eq{refinedtor}
    \rho_{{}_\Gam} \ = \ \rho_{{}_{C^\b,\Gam}} \ := \  \phi_{C^\b}(c_{{}_\Gam}) \ \in \Det\big(\,H^\b(\pa)\,\big),
\end{equation}
where $ \phi_{C^\b}$ is the canonical map \refe{isomorphism}.
\edefe

\subsection{The odd signature operator}\Label{SS:oddsign}
Let $M$  be a smooth closed oriented manifold of odd dimension $d=2r-1$ and let $(E,\n)$ be a flat vector bundle over $M$. We denote by $\Ome^k(M,E)$
the space of smooth differential forms on $M$  of degree $k$ with values in $E$ and by
\[
    \n:\, \Ome^\b(M,E) \ \longrightarrow \Ome^{\b+1}(M,E)
\]
the covariant differential induced by the flat connection on $E$.

Fix a Riemannian metric $g^M$ on $M$ and let $*:\Ome^\b(M,E)\to \Ome^{d-\b}(M,E)$ denote the Hodge $*$-operator. Define the {\em chirality operator}
$\Gam= \Gam(g^M):\Ome^\b(M,E)\to \Ome^\b(M,E)$ by the formula
\eq{Gam}
    \Gam\, \ome \ := \ i^r\,(-1)^{\frac{k(k+1)}2}\,*\,\ome, \qquad \ome\in \Ome^k(M,E),
\end{equation}
with $r$ given as above by \/ $r=\frac{d+1}2$. The numerical factor in \refe{Gam} has been chosen so that $\Gam^2=1$, cf. Proposition~3.58 of
\cite{BeGeVe}.

\defe{oddsign}
The {\em odd signature operator} is the operator
\begin{equation} \Label{E:oddsignGam}
    \B\ = \  \B(\n,g^M) \ := \ \Gam\,\n \ + \ \n\,\Gam:\,\Ome^\b(M,E)\ \longrightarrow \  \Ome^\b(M,E).
\end{equation}
We denote by $\B_k$ the restriction of\/ $\B$ to the space $\Ome^{k}(M,E)$.
\edefe

\subsection{The graded determinant of the odd signature operator}\label{SS:grdet}
Note that for each $k=0\nek d$, the operator $\B^2$ maps $\Ome^k(M,E)$ into itself. Suppose $\calI$ is an interval of the form $[0,\lam],\
(\lam,\mu]$, or $(\lam,\infty)$ ($\mu>\lam\ge0$). Denote by $\Pi_{\B^2,\calI}$ the spectral projection of $\B^2$ corresponding to the set of
eigenvalues, whose absolute values lie in $\calI$. Set
\eq{OmecalI}\notag
    \Ome^\b_{\calI}(M,E)\ := \ \Pi_{\B^2,\calI}\big(\, \Ome^\b(M,E)\,\big)\ \subset\  \Ome^\b(M,E).
\end{equation}
If the interval $\calI$ is bounded, then, cf. Section~6.10 of \cite{BrKappelerRATdetline}, the space $\Ome^\b_\calI(M,E)$ is finite dimensional.

For each $k=0\nek d$, set
\eq{ome+-}
 \begin{aligned}
  \Ome^k_{+,\calI}(M,E) \ &:= \ \Ker\,(\n\,\Gam)\,\cap\,\Ome^k_{\calI}(M,E) \ = \ \big(\,\Gam\,(\Ker\,\n)\,\big)\,\cap\,\Ome^k_{\calI}(M,E);\\
  \Ome^k_{-,\calI}(M,E) \ &:= \ \Ker\,(\Gam\,\n)\,\cap\,\Ome^k_{\calI}(M,E)
  \ = \ \Ker\,\n\,\cap\,\Ome^k_{\calI}(M,E).
 \end{aligned}
\end{equation}
Then
\eq{Ome>0=directsum}
    \Ome^k_{\calI}(M,E) \ = \ \Ome^k_{+,\calI}(M,E) \,\oplus \Ome^k_{-,\calI}(M,E) \qquad\text{if}\quad 0\not\in \calI.
\end{equation}
We consider the decomposition \refe{Ome>0=directsum} as a {\em grading} \footnote{Note, that our grading is opposite to the one considered in
\cite[\S2]{BFK3}.} of the space $\Ome^\b_{\calI}(M,E)$, and refer to $\Ome^k_{+,\calI}(M,E)$ and $\Ome^k_{-,\calI}(M,E)$ as the positive and negative
subspaces of\/ $\Ome^k_{\calI}(M,E)$.

Set
\[
    \Ome^{\even}_{\pm,\calI}(M,E)\ =\ \bigoplus_{p=0}^{r-1}\, \Ome^{2p}_{\pm,\calI}(M,E)
\]
and let $\B^\calI$ and $\B^\calI_\even$ denote the restrictions of $\B$ to the subspaces $\Ome^\b_{\calI}(M,E)$ and $\Ome^\even_{\calI}(M,E)$
respectively. Then $\B^\calI_{\even}$ maps $\Ome^{\even}_{\pm,\calI}(M,E)$ to itself. Let $\B_{{\even}}^{\pm,\calI}$ denote the restriction of
$\B^\calI_{{\even}}$ to the space $\Ome^{\even}_{\pm,\calI}(M,E)$. Clearly, the operators $\B_{{\even}}^{\pm,\calI}$ are bijective whenever $0\not\in
\calI$.

\defe{grdetB}
Suppose $0\not\in \calI$. The {\em graded determinant} of the operator $\B^{\calI}_{\even}$ is defined by
\eq{grdetB}
    \Detgrtetnp(\B^{\calI}_{\even}) \ := \ \frac{ \Det_\tet(\B^{+,\calI}_{\even})}{\Det_\tet(-\B^{-,\calI}_{\even})}
    \ \ \in\ \ \CC\backslash \{0\},
\end{equation}
where $\Det_\tet$ denotes the $\zet$-regularized determinant associated to the Agmon angle $\tet\in (-\pi,0)$, cf., for example, \S6 of\/
\cite{BrKappelerRATdetline}.
\edefe
It follows from formula (6.17) of \cite{BrKappelerRATdetline} that \refe{grdetB} is independent of the choice of $\tet\in (-\pi,0)$.

\subsection{The refined analytic torsion}\Label{SS:rho}
Since the covariant differentiation $\n$ commutes with $\B$, the subspace $\Ome^\b_\calI(M,E)$ is a subcomplex of the twisted de Rham complex
$(\Ome^\b(M,E),\n)$. Clearly, for each $\lam\ge0$, the complex $\Ome^\b_{(\lam,\infty)}(M,E)$ is acyclic. Since
\eq{Ome=Ome0+Ome>0}
    \Ome^\b(M,E) \ = \  \Ome^\b_{[0,\lam]}(M,E)\,\oplus\,\Ome^\b_{(\lam,\infty)}(M,E),
\end{equation}
the cohomology $H^\b_{[0,\lam]}(M,E)$ of the complex $\Ome^\b_{[0,\lam]}(M,E)$ is naturally isomorphic to the cohomology $H^\b(M,E)$.

Let $\Gam_{\hskip-1pt{}\calI}$ denote the restriction of $\Gam$ to $\Ome^\b_\calI(M,E)$. For each $\lam\ge0$, let
\eq{thoGam0lam}
    \rho_{{}_{\Gam_{\hskip-1pt{}_{[0,\lam]}}}}\ = \ \rho_{{}_{\Gam_{\hskip-1pt{}_{[0,\lam]}}}}(\n,g^M) \in \ \Det(H^\b_{[0,\lam]}(M,E))
\end{equation}
denote the refined torsion of the finite dimensional complex $(\Ome^\b_{[0,\lam]}(M,E),\n)$ corresponding to the chirality operator
$\Gam_{\hskip-1pt{}_{[0,\lam]}}$, cf. \refd{refinedtorsion}. We view $\rho_{{}_{\Gam_{\hskip-1pt{}_{[0,\lam]}}}}$ as an element of $\Det(H^\b(M,E))$
via the canonical isomorphism between $H^\b_{[0,\lam]}(M,E)$ and $H^\b(M,E)$.

It is shown in Proposition~7.8 of \cite{BrKappelerRATdetline} that the nonzero element
\eq{rho}
    \rho(\n) \ = \ \rho(\n,g^M) \ := \ \Detgrtetnp(\B^{(\lam,\infty)}_{\even})\cdot \rho_{{}_{\Gam_{\hskip-1pt{}_{[0,\lam]}}}}
    \ \in  \ \Det(H^\b(M,E))
\end{equation}
is independent of the choice of $\lam\ge0$. Further, $\rho(\n)$ is independent of the choice of the Agmon angle $\tet\in (-\pi,0)$ of\/ $\B_\even$.

If the odd signature operator is invertible then $\Det(H^\b(M,E))$ is canonically isomorphic to $\CC$ and $\rho_{{}_{\Gam_{\hskip-1pt{}_{\{0\}}}}}=
1$. Hence, $\rho(\n)$ is a complex number which coincides with the graded determinant $\Det_{\operatorname{gr},\tet}(\B_{\even})=
\Det_{\operatorname{gr},\tet}(\B^{(0,\infty)}_{\even})$. This case was studied in \cite{BrKappelerRAT}.

Let $\B_\trivial(g^M):\Ome^\even(M,\CC)\to \Ome^\even(M,\CC)$ denote the even part of the odd signature operator
$\Gam{}d+d\Gam:\Ome^\b(M,\CC)\to \Ome^\b(M,\CC)$ corresponding to the trivial line bundle $M\times\CC\to M$.

\defe{rat}
The {\em refined analytic torsion} is the element
\eq{metricindep}
    \rat(\n)\  := \ \rho(\n,g^M)\cdot e^{i\pi\cdot\rank E\cdot\eta_\trivial(g^M)}\ \in \ \Det(H^\b(M,E)),
\end{equation}
where $g^M$ is any Riemannian metric on $M$, $\rho(\n,g^M)\in \Det(H^\b(M,E))$ is defined in \refe{rho}, and
\[
    \eta_\trivial(g^M) \ = \ \frac12\,\eta(0,\B_\trivial),
\]
is one half of the value at zero of the $\eta$-function of the operator $\B_\trivial$.

In particular, if $\dim M\equiv 1\ \ (\MOD 4)$, then $\eta_\trivial(g^M)=0$, cf. \cite{APS1}, and $\rat(\n)= \rho(\n,g^M)$.
\edefe
It is shown in Theorem~9.6 of \cite{BrKappelerRATdetline} that $\rat(\n)$ is independent of $g^M$.

\section{The Determinant Line Bundle over the Space of Representations}\label{S:detlinebundle}

The space $\Rep$ of complex $n$-dimensional representations of $\p$ has a natural structure of a complex analytic space, cf., for example,
\cite[\S13.6]{BrKappelerRAT}. The disjoint union
\eq{Detbundle}
    \Dbundle \ := \ \bigsqcup_{\alp\in\Rep}\, \Det\big(H^\b(M,E)\big)
\end{equation}
has a natural structure of a holomorphic line bundle over $\Rep$, called the {\em determinant line bundle}. In this section we describe this
structure, using a CW-decomposition of $M$. Then, by construction, the Farber-Turaev torsion $\rtt(\alp)$ is a nowhere vanishing holomorphic
section of $\Dbundle$. In particular, it defines a holomorphic trivialization of $\Dbundle$. Note, however, that this trivialization is not
canonical since it depends on the Euler structure $\eps$.

\subsection{The flat vector bundle induced by a representation}\label{SS:Ealp}
Denote by $\pi:\tilM\to M$ the universal cover of $M$ and by $\pi_1(M)$ the fundamental group of $M$, viewed as the group of deck transformations of
$\tilM\to M$. For $\alp\in \Rep$, we denote by
\eq{Ealp}
  E_\alp \ := \ \tilM\times_\alp \CC^n \ \longrightarrow M
\end{equation}
the flat vector bundle induced by $\alp$. Let $\n_\alp$ be the flat connection on $E_\alp$ induced from the trivial connection on $\tilM\times\CC^n$.
We will also denote by $\n_\alp$ the induced differential
\[
    \n_\alp:\, \Ome^\b(M,E_\alp) \ \longrightarrow \Ome^{\b+1}(M,E_\alp),
\]
where $\Ome^\b(M,E_\alp)$ denotes the space of smooth differential forms of $M$ with values in $E_\alp$.

For each connected component (in classical topology) $\calC$ of $\Rep$, all the bundles $E_\alp$, $\alp\in \calC$, are isomorphic, see e.g.
\cite{GoldmanMillson88}.

\subsection{The combinatorial cochain complex}\label{SS:CKalp}
Fix a CW-decomposition $K=\{e_1\nek e_N\}$ of $M$. For each $j=1\nek N$, fix a lift $\tile_j$, i.e., a cell of the CW-decomposition of $\tilM$, such
that $\pi(\tile_j)= e_j$. By \refe{Ealp}, the pull-back of the bundle $E_\alp$ to $\tilM$ is the trivial bundle $\tilM\times\CC^n\to \tilM$. Hence,
the choice of the cells $\tile_1\nek \tile_N$ identifies the cochain complex $C^\b(K,\alp)$ of the CW-complex $K$ with coefficients in $\Ea$ with the
complex
\eq{C(K,alp)}
    \begin{CD}
        0 \ \to \CC^{n\cdot k_0} @>{\pa_0(\alp)}>> \CC^{n\cdot k_1} @>{\pa_1(\alp)}>>\cdots @>{\pa_{d-1}(\alp)}>> \CC^{n\cdot k_d} \ \to \ 0,
    \end{CD}
\end{equation}
where $k_j\in \ZZ_{\ge0}$  ($j=0\nek d$) is equal to the number of $j$-dimensional cells of $K$ and the differentials $\pa_j(\alp)$ are
$(nk_{j}\times{}nk_{j-1})$-matrices depending analytically on $\alp\in \Rep$.

The cohomology of the complex \refe{C(K,alp)} is canonically isomorphic to $H^\b(M,\Ea)$. Let
\eq{phialp}
    \phi_{C^\b(K,\alp)}:\, \Det\big(\,C^\b(K,\alp)\,\big) \ {\longrightarrow} \ \Det\big(\,H^\b(M,\Ea)\,\big)
\end{equation}
denote the isomorphism \refe{isomorphism}.

\subsection{The holomorphic structure on $\Dbundle$}\label{SS:holonDet}
The standard bases of $\CC^{n\cdot k_j}$ ($j=0\nek d$) define an element $c\in \Det\big(C^\b(K,\alp)\big)$, and, hence, an isomorphism
\[
    \psi_\alp:\,\CC \ {\longrightarrow} \ \Det\big(C^\b(K,\alp)\big), \qquad z\ \mapsto z\cdot c.
\]
Then the map
\eq{combsection}
    \sig:\, \alp \ \mapsto \ \phi_{C^\b(K,\alp)}\big(\,\psi_\alp(1)\,\big) \ \in \ \Det\big(\,H^\b(M,\Ea)\,\big), \qquad \alp\in \Rep
\end{equation}
is a nowhere vanishing section of the determinant line bundle $\Dbundle$ over $\Rep$.
\defe{holsection}
We say that a section $s(\alp)$ of\/ $\Dbundle$ is {\em holomorphic} if there exists a holomorphic function $f(\alp)$ on $\Rep$, such that $s(\alp)=
f(\alp)\cdot{}\sig(\alp)$.
\edefe
This defines a holomorphic structure on $\Dbundle$, which is independent of the choice of the lifts $\tile_1\nek \tile_N$ of $e_1\nek e_N$, since for
a different choice of lifts the section $\sig(\alp)$ will be multiplied by a constant. In the next subsection we show that this holomorphic structure
is also independent of the CW-decomposition $K$ of $M$.

\subsection{The Farber-Turaev torsion}\label{SS:Turaevtorsion}
The choice of the lifts $\tile_1\nek \tile_N$ of $e_1\nek e_N$ determines an {\em Euler structure} on $M$, while the ordering of the cells
$e_1\nek e_N$ determines a cohomological orientation $\gro$, cf. \cite[\S20]{Turaev01}. Moreover, every Euler structure and every cohomological
orientation can be obtained in this way. The Farber-Turaev torsion $\rtt(\alp)$, corresponding to the pair $(\eps,\gro)$, is, by definition,
\cite[\S6]{FarberTuraev00}, equal to the element $\sig(\alp)$ defined in \refe{combsection}. In particular, it is a non-vanishing holomorphic
section of $\Dbundle$, according to \refd{holsection}. Since the Farber-Turaev torsion is independent of the choice of the CW-decomposition of
$M$ \cite{Turaev90,FarberTuraev00}, so is the holomorphic structure defined in \refd{holsection}.

\subsection{The acyclic case}\label{SS:Turaevacyclic}
If the representation $\alp$ is acyclic, i.e., $H^\b(M,E_\alp)=0$, then the determinant line $\Det(H^\b(M,E_\alp))$ is canonically isomorphic to
$\CC$. Hence, the Farber-Turaev torsion can be viewed as a complex-valued function on the set $\Repo\subset \Rep$ of acyclic representations. It
is easy to see, cf. Theorem~4.3 of \cite{BurgheleaHaller_Euler}, that this function is holomorphic on $\Repo$. (Moreover, it is a rational
function on $\Rep$, all whose poles are in $\Rep\backslash\Repo$). In particular, the holomorphic structure on $\Dbundle$, which we defined
above, coincides, when restricted to $\Repo$, with the natural holomorphic structure obtained from the canonical isomorphism
$\Dbundle|_{\Repo}\simeq \Repo\times\CC$.

\smallskip
We summarize the results of this section in the following
\prop{combholstr}
a. \ The holomorphic structure defined in \refd{holsection} is independent of any choices made.

b. \ For every Euler structure $\eps$ and every cohomological orientation $\gro$, the Farber-Turaev torsion $\rtt(\alp)$ is a holomorphic
section of the determinant line bundle $\Dbundle$.

c. \ The restriction of $\rtt(\alp)$ to the open subset $\Repo\subset \Rep$ of acyclic representations is a holomorphic function.
\eprop

\section{Refined Analytic Torsion as a Holomorphic Section}\label{S:holomorphic}

One of the main results of this paper is that the refined analytic torsion $\rat$ is a non-vanishing holomorphic section of\/ $\Dbundle$. More
precisely, the following theorem holds.
\th{holomsection}
The refined analytic torsion $\rat$ is a holomorphic section of the determinant bundle $\Dbundle$, i.e., for any Euler structure $\eps$ and any
cohomological orientation $\gro$, the ratio $\rat/\rtt$ is a holomorphic function on $\Rep$.

In particular, the restriction of $\rat$ to the set $\Repo$ of acyclic representations, viewed as a complex-valued function via the canonical
isomorphism
\[
    \Dbundle|_{\Repo}\ \simeq\ \Repo\times\CC,
\]
is a holomorphic function on $\Repo$.
\eth

We prove this theorem in two steps: in this section we show that $\rat$ is holomorphic on $\Rep\backslash\Sig(M)$, where $\Sig(M)$ is the set of
singular points of the complex analytic set $\Rep$. In the next section we will use this result to calculate the ratio of the refined analytic
and the Farber-Turaev torsions. This calculation and the fact that the Farber-Turaev torsion is holomorphic, will imply that $\rat$ is
holomorphic everywhere, cf. \refss{prholomsection}.

The main result of this section is the following
\prop{holnearalp}
Let $\alp_0\in \Rep\backslash\Sig(M)$. Then the refined analytic torsion $\rat$, viewed as a section of\/  $\Dbundle$, is holomorphic in a
neighborhood of $\alp_0$ with respect to the holomorphic structure defined in \refs{detlinebundle}.
\eprop

For convenience of the reader and in order to illustrate the main ideas of the proof we, first, prove the proposition for the case when $\alp_0$ is
acyclic. Then in a neighborhood of $\alp_0$ the refined analytic torsion can be viewed as a complex-valued function, and we shall show that this
function is holomorphic at $\alp_0$.

\subsection{Reduction to a finite dimensional complex}\label{SS:redtofd}
Fix a Riemannian metric $g^M$ on $M$ and a number $\lam\ge0$ such that there are no eigenvalues of\/ $\B(\n_{\alp_0},g^M)^2$ with absolute value
equal to $\lam$. Then there exists a neighborhood $U_\lam\subset \Repo\backslash\Sig(M)$ of $\alp_0$ such that the same property holds for all
$\alp\in U_\lam$. By Proposition~13.2 of \cite{BrKappelerRAT} the function $\alp\mapsto \Detgrtetnp(\B_\even^{(\lam,\infty)}(\na,g^M))$ is
holomorphic\footnote{Proposition~13.2 of \cite{BrKappelerRAT} only deals with the case where $\B$ is invertible and $\lam=0$. But a verbatim
repetition of the same proof with $\B$ replaced everywhere by $\B^{(\lam,\infty)}$ works in our more general situation.} on $U_\lam$. It follows
now from \refe{rho} and \refe{metricindep} that to prove \refp{holnearalp} it is enough to show that the function
\[
    \alp \mapsto \ \rho_{{}_{\Gam_{\hskip-1pt{}_{[0,\lam]}}}}(\na)\ =\ \rho_{{}_{\Gam_{\hskip-1pt{}_{[0,\lam]}}}}(\na,g^M)
\]
is holomorphic.

\subsection{Reduction to one-parameter families of representations}\label{SS:oneparameter}
By Hartog's theorem, \cite[Th.~2.2.8]{HormanderSCV}, it is enough to show that for every holomorphic curve $\gam:\calO\to U_\lam$, where $\calO$ is a
connected open neighborhood of 0 in $\CC$, such that $\gam(0)= \alp_0$,
\[
    z \ \mapsto \rho_{{}_{\Gam_{\hskip-1pt{}_{[0,\lam]}}}}(\n_{\gam(z)}), \qquad z\in \calO,
\]
is a holomorphic function on $\calO$.

\subsection{A family of connections}\label{SS:familycon}
Let us introduce some additional notations. Let $E$ be a vector bundle over $M$ and let $\n$ be a flat connection on $E$. Fix a base point $x_*\in M$
and let $E_{x_*}$ denote the fiber of $E$ over $x_*$. We will identify $E_{x_*}$ with $\CC^n$ and $\pi_1(M,x_*)$ with $\pi_1(M)$.

For a closed  path $p:[0,1]\to M$ with $p(0)=p(1)= x_*$, we denote by $\Mon_\n(p)\in \End{}E_{x_*}\simeq \Mat_{n\times n}(\CC)$ the monodromy of $\n$
along $p$. Since the connection $\n$ is flat, $\Mon_\n(p)$ depends only on the class $[p]$ of $p$ in $\pi_1(M)$. Hence, the map $p\mapsto
\Mon_\n(\phi)$ defines an element of $\Rep$, called  the {\em monodromy representation} of $\n$.

Suppose now that $\calO\subset \CC$ is a connected open set. Let $\gam:\calO\to \Rep$ be a holomorphic curve. By Proposition~4.5 of
\cite{GoldmanMillson88}, all the bundles $E_{\gam(z)}, \ z\in \calO$, are isomorphic to each other. In other words, there exists a vector bundle
$E\to M$ and a family of flat connections $\n_z,\ z\in \calO$, on $E$, such that the monodromy representation of $\n_z$ is equal to $\gam(z)$ for all
$z\in \calO$. Moreover, Lemma~B.6 of \cite{BrKappelerRAT} shows that, for any $z_0\in \calO$, the family $\n_z$ can be chosen so that there exists a
one-form $\ome\in \Ome^1(M,\End{E})$ such that
\eq{nz-nz0}
    \n_z \ = \ \n_{z_0} \ + \ (z-z_0)\,\ome \ + \ o(z-z_0).
\end{equation}
Since $z_0\in \calO$ is arbitrary, it follows now from the discussion in \refss{oneparameter}, that to finish the proof of \refp{holnearalp} we only
need to show that the function
\eq{f(z)}
    f(z) \ := \ \rho_{{}_{\Gam_{\hskip-1pt{}_{[0,\lam]}}}}(\n_z,g^M)
\end{equation}
is complex differentiable at $z_0$, i.e., there exists $a\in \CC$, such that
\[
    f(z)\ =\ f(z_0)\ +\ (z-z_0)\,a\ +\ o(z-z_0).
\]

\subsection{Choice of a basis}\label{SS:basis}
Let $\Pi_{[0,\lam]}(z)$ ($z\in \calO$) denote the spectral projection of the operator $\B(\n_z,g^M)^2$, corresponding to the set of eigenvalues of
$\B(\n_z,g^M)^2$, whose absolute value is $\le\lam$, cf. \refss{redtofd}. It follows from \refe{nz-nz0} that there exists a bounded operator
$P:\Ome^\b(M,E)\to \Ome^\b(M,E)$ such that
\eq{Piz-Piz0}
    \Pi_{[0,\lam]}(z) \ = \ \Pi_{[0,\lam]}(z_0) \ + \ (z-z_0)\, P \ +\ o(z-z_0).
\end{equation}

We denote by $\Ome^\b(z)$ the image of $\Pi_{[0,\lam]}(z)$. For each $j=0\nek r-1$, fix a basis $\bfw_j= \{w^1_j\nek w^{m_j}_j\}$ of\/ $\Ome^j(z_0)$
and set $\bfw_{d-j}:= \{\Gam{}w^1_j\nek \Gam{}w^{m_j}_j\}$. To simplify the notation we will write $\bfw_{d-j}= \Gam\bfw_j$. Then $\bfw_j$ is a basis
for $\Ome^j(z_0)$ for all $j=0\nek d$.

For each $z\in \calO$, $j=0\nek d$, set
\[
    \bfw_j(z)\ = \ \big\{\,w^1_j(z)\nek w^{m_j}_j(z)\,\big\} \ := \ \big\{\,\Pi_{[0,\lam]}(z)\,w^1_j\nek \Pi_{[0,\lam]}(z)\,w^{m_j}_j\,\big\}.
\]
It follows from the definition of $U_\lam$ that the projection $\Pi_{[0,\lam]}(z)$ depends continuously on $z$. Hence, there exists a neighborhood
$\calO'\subset \calO$ of $z_0$, such that $\bfw_j(z)$ is a basis of $\Ome^j(z)$ for all $z\in \calO'$, $j=0\nek d$. Further, since
$\Pi_{[0,\lam]}(z)$ commutes with $\Gam$, we obtain
\eq{oome=Gamoome}
    \bfw_{d-j}(z)\ =\ \Gam\,\bfw_j(z).
\end{equation}
Clearly, $\bfw_j(z_0)= \bfw_j$ for all $j=0\nek d$.

Let
\[
    \phi_{\Ome^\b(z)}:\, \Det\big(\,\Ome^\b(z)\,\big) \ \longrightarrow \ \Det\big(\,H^\b(M,E_{\gam(z)})\,\big) \ \simeq \ \CC
\]
denote the isomorphism \refe{isomorphism}. For $z\in \calO'$, let $w(z)\in \Det\big(\,\Ome^\b(z)\,\big)$  be the element determined by the basis
$\bfw_1(z)\nek \bfw_d(z)$ of $\Ome^\b(z)$. More precisely, we introduce
\[
    w_j(z) \ =\  w_j^1(z)\wedge\cdots\wedge{}w_j^{m_j}(z) \ \in \ \Det\big(\,\Ome^j(z)\,\big),
\]
and set
\[
    w(z) \ := \ w_0(z)\otimes{}w_1(z)^{-1}\otimes\cdots\otimes{}w_d(z)^{-1}.
\]
Then, according to \refd{refinedtorsion}, it follows from \refe{oome=Gamoome} that, for all $z\in \calO'$, the refined torsion of the complex
$\Ome^\b(z)$ is equal to $\phi_{\Ome^\b(z)}(w(z))$, i.e.,
\eq{rho0lam=}
    \rho_{{}_{\Gam_{\hskip-1pt{}_{[0,\lam]}}}}(\n_z) \ = \ \phi_{\Ome^\b(z)}(w(z)).
\end{equation}

\subsection{Reduction to a family of differentials}\label{SS:familyofd}
For each $z\in \calO'$, the space $\Ome^\b(z)$ is a subcomplex of $\big(\Ome^\b(M,E),\n_z\big)$, whose cohomology is canonically isomorphic to the
cohomology of $\big(\Ome^\b(M,E),\n_z\big)$ and, hence, to $H^\b(M,E_{\gam(z)})$. Using the basis $\bfw_j(z)$ we define the isomorphism
\[
    \psi_j(z):\, \CC^{m_j}  \ {\longrightarrow}\ \Ome^j_{[0,\lam]}(z)
\]
by the formula
\eq{psi(z)}
    \psi_j(z)(x_1\nek x_{m_j}) \ := \ \sum_{k=1}^{m_j}\,x_j\,w_j^k(z) \ = \ \sum_{k=1}^{m_j}\,x_j\,\Pi_{[0,\lam]}(z)\,w_j^k.
\end{equation}
We conclude that for each $z\in \calO'$, the complex $\big(\Ome^\b(z),\n_z\big)$ is isomorphic to the complex
\eq{coordOme0lam}
    \begin{CD}
       (W^\b,d(z)):\qquad 0 \ \to \CC^{m_0} @>{d_0(z)}>> \CC^{m_1} @>{d_1(z)}>>\cdots @>{d_{d-1}(z)}>> \CC^{m_d} \ \to \ 0,
    \end{CD}
\end{equation}
where
\eq{d(z)}
    d_j(z) \ := \ \psi_{j+1}(z)^{-1}\circ \n_z\circ \psi_j(z), \qquad j=0\nek d.
\end{equation}
It follows from \refe{Piz-Piz0} and \refe{psi(z)} that $d_j(z)$ is complex differentiable at $z_0$, i.e., there exists a
$(m_{j+1}\times{}m_j)$-matrix $A$ such that
\[
    d_j(z) \ = \ d_j(z_0) \ + \ (z-z_0)\,A \ + \ o(z-z_0).
\]
Let $\psi(z):= \bigoplus_{j=0}^d\psi_j(z)$. Since $\Gam\big(\Ome^j(z)\big)= \Ome^{d-j}(z)$ ($j=0\nek d$), we conclude that $m_j= m_{d-j}$. From
\refe{oome=Gamoome} we obtain that
\eq{tilGam}
    \tilGam\ :=\  \psi^{-1}(z)\circ \Gam \circ \psi(z)
\end{equation}
is independent of $z\in \calO'$ and
\eq{Gamofbasis}
    \tilGam:\, (x_1\nek x_{m_j}) \ \mapsto (x_1\nek x_{m_j}), \qquad j=0\nek d.
\end{equation}
It follows from \refe{d(z)} and \refe{tilGam} that
\eq{Gam-tilGam}
    \rho_{{}_{\tilGam}}(z) \ = \ \rho_{{}_{\Gam_{\hskip-1pt{}_{[0,\lam]}}}}(\n_z),
\end{equation}
where $\rho_{{}_{\tilGam}}(z)$ denotes the refined torsion of the finite dimensional complex $(W^\b,d(z))$ corresponding to the chirality operator
$\tilGam$, cf. \refd{refinedtorsion}.

Let $\phi_{W^\b}(z):\Det(W^\b)\to \Det\big(H^\b(d(z))\big)$ denote the isomorphism \refe{isomorphism}. The standard bases of $\CC^{m_j}$ ($j=0\nek
d$) define an element $\tilw\in \Det(W^\b)$. From \refe{Gamofbasis}  and the definition \refe{refinedtor} of $\rho_{\tilGam}(z)$ we conclude that
\eq{rhotilgam}
    \rho_{{}_{\tilGam}}(z) \ = \ \phi_{W^\b}(z)(\tilw).
\end{equation}

Hence, to finish the proof of \refp{holnearalp} in the case when $\alp_0$ is acyclic it remains to show that the function $z\mapsto
\phi_{W^\b}(z)(\tilw)$ is complex differentiable at $z_0$. In view of \refe{d(z)}, this follows from the following
\lem{asyclicanalytic}
Let
\eq{Canalytic}
    \begin{CD}
       (C^\b,\pa(z)):\qquad 0 \ \to \CC^{n\cdot k_0} @>{\pa_0(z)}>> \CC^{n\cdot k_1}
       @>{\pa_1(z)}>>\cdots @>{\pa_{d-1}(z)}>> \CC^{n\cdot k_d} \ \to \ 0,
    \end{CD}
\end{equation}
be a family of acyclic complexes defined for all $z$ in an open set $\calO\subset \CC$. Suppose that the differentials $\pa_j(z)$ are complex
differentiable at $z_0\in \calO$. Then for any $c\in \Det(C^\b)$ the function $z\mapsto \phi_{(C^\b,\pa(z))}(c)$ is complex differentiable at $z_0$.
\elem
\prf
It is enough to prove the lemma for one particular choice of $c$. To make such  a choice let us fix for each $j=0\nek d$ a complement of
$\IM(\pa_{j-1}(z_0))$ in $C^j$ and a basis $v_j^1\nek v_j^{l_j}$ of this complement. Since the complex $C^\b$ is acyclic, for all $j=0\nek d$, the
vectors
\eq{basisvz0}
    \pa_{j-1}(z_0)\,v_{j-1}^1\,,\ \ldots\,,\  \pa_{j-1}(z_0)\,v_{j-1}^{l_{j-1}},\ \,  v_j^1\ ,\,\ldots\,,\  v_j^{l_j}
\end{equation}
form a basis of $C^j$. Let $c\in \Det(C^\b)$ be the element defined by these bases. Then, for all $z$ close enough to $z_0$  and for all $j=0\nek d$,
\eq{basisvz}
    \pa_{j-1}(z)\,v_{j-1}^1\,,\ \ldots\,,\  \pa_{j-1}(z)\,v_{j-1}^{l_{j-1}},\ \, v_j^1\,,\ \ldots\,,\   v_j^{l_j}
\end{equation}
is also a basis of $C^j$. Let $A_j(z)\ (j=0\nek d)$ denote the non-degenerate matrix transforming the basis \refe{basisvz} to the basis
\refe{basisvz0}. Then, by the definition of the isomorphism $\phi_{(C^\b,\pa(z))}$, cf. \S2.4 of \cite{BrKappelerRATdetline},
\eq{phi(z)(c)}
    \phi_{(C^\b,\pa(z))}(c)\ = \ (-1)^{\calN(C^\b)}\ \prod_{j=0}^d\,\Det(A(z))^{(-1)^j},
\end{equation}
where $\calN(C^\b)$ is the integer defined in formula (2.15) of \cite{BrKappelerRATdetline} which is independent of $z$. Clearly, the matrix valued
functions $A_j(z)$ and, hence, their determinants are complex differentiable at $z_0$. Thus, so is the function $z\mapsto \phi_{(C^\b,\pa(z))}(c)$.
\eprf

\subsection{Sketch of the proof of Proposition~4.2 in the non-acyclic case}\label{SS:sketchnonacyclic}
Let $\n_z$ be the family of connections \refe{nz-nz0}. To prove \refp{holnearalp} in the case when $\alp_0$ is not acyclic it is enough to show
that the function
\eq{f(z)2}
    f(z) \ := \ \frac{\rho_{{}_{\Gam_{\hskip-1pt{}_{[0,\lam]}}}}(\n_z,g^M)}{\rtt(\gam(z))}
\end{equation}
is complex differentiable at $z_0$. Here $\rtt(\gam(z))$ stands for the Farber-Turaev torsion associated to the representation $\gam(z)$, the Euler structure $\eps$
and the cohomological orientation $\gro$, cf. \refss{Turaevtorsion}. To see this we consider the integration map
\[
    J_z:\, \Ome^\b(z)\subset \Ome^\b(M,E) \ \longrightarrow \ C^\b(K,\gam(z)),
\]
where $C^\b(K,\gam(z))$ is the cochain complex corresponding to the CW-decomposition $K= \linebreak\{e_1\nek e_N\}$, cf. \refss{CKalp}. Note that the
integration of $E$-valued differential forms is defined using a trivialization of $E$ over each cell $e_j$, and, hence, it depends on the flat
connection $\n_z$, cf. below. We then consider the cone complex $\Cone^\b(J_z)$ of the map $J_z$. This is a finite dimensional acyclic complex with a
fixed basis, obtained from the bases of $\Ome^\b(z)$ and $C^\b(K,\gam(z))$. The torsion of this complex is equal to $f(z)$. An application of
\refl{asyclicanalytic} to this complex proves \refp{holnearalp}.

In the definition of the integration map $J_z$ we have to take into account the fact that the vector bundles $E_{\gam(z)}$ and $E$ are isomorphic but
not equal. The standard integration map, cf. \refss{combin}, is a map from $\Ome^\b(M,E)$ to the cochain complex $C^\b(K,E)$ of $K$ with coefficients
in $E$, which is not equal to the complex $C^\b(M,E_{\gam(z)})$. There is a natural isomorphism between the complexes $C^\b(K,E)$ and
$C^\b(K,\gam(z))$ which depends on $z$. The study of this isomorphism, which is conducted in \refss{tilC-C}, is important for the understanding of
the properties of $J_z$. In particular, it is used to show that $J_z$ is complex differentiable at $z_0$, which implies that the cone complex
$\Cone^\b(J_z)$ satisfies the conditions of \refl{asyclicanalytic}.

\subsection{The cochain complex of the bundle $E$}\label{SS:combin}
Fix a CW-decomposition $K= \{e_1\nek e_N\}$ of $M$. For each $j=1\nek N$ choose a point $x_j\in e_j$ and let $E_{x_j}$ denote the fiber of $E$ over
$x_j$. The cochain complex  of the CW-decomposition $K$ with coefficients in the flat bundle $(E,\n_z)$ can be identified with the complex
$(C^\b(K,E),\pa^\prime(z))$
\eq{tilC(K,alp)}
    \begin{CD}\displaystyle
       0 \ \to \bigoplus_{\dim e_i=0} E_{x_i} @>{\pa^\prime_0(z)}>> \displaystyle\bigoplus_{\dim e_i=1} E_{x_i}
       @>{\pa^\prime_1(z)}>>\cdots @>{\pa^\prime_{d-1}(z)}>> \displaystyle\bigoplus_{\dim e_i=d} E_{x_i} \ \to \ 0.
    \end{CD}
\end{equation}
We use the prime in the notation of the differentials $\pa^\prime_j$ in order to distinguish them from the differentials of the cochain complex
$C^\b(K,\gam(z))$ defined in \refe{C(K,alp)}.

It follows from \refe{nz-nz0} that $\pa^\prime_j(z)$ are complex differentiable at $z_0$, i.e., there exist linear maps
\[
    a_j:\,\bigoplus_{\dim e_i=j}\,E_{x_i}\ \longrightarrow\ \bigoplus_{\dim e_i=j+1}\,E_{x_i}
\]
such that
\eq{tilpa'(z)}\notag
    \pa^\prime_j(z) \ = \ \pa^\prime_j(z_0)\ + \ (z-z_0)\,a_j\ + \ o(z-z_0), \qquad j=1\nek d-1.
\end{equation}

\subsection{Relationship with the complex $C^\b(K,\gam(z))$}\label{SS:tilC-C}
Recall that for each $z\in \calO'$ the monodromy representation of $\n_z$ is equal to $\gam(z)$. Let $\pi:\tilM\to M$ denote the universal cover of
$M$ and let $\tilE= \pi^*E$ denote the pull-back of the bundle $E$ to $\tilM$. Recall that in \refss{familycon} we fixed a point $x_*\in M$. Let
$\tilx_*\in \tilM$ be a lift of $x_*$ to $\tilM$ and fix a basis of the fiber $\tilE_{\tilx_*}$ of $\tilE$ over $x_*$. Then, for each $z\in \calO'$,
the flat connection $\n_z$ identifies $\tilE$ with the product $\tilM\times \CC^n$. Let $\tile_j \ (j=1\nek N)$ be the lift of the cell $e_j$ fixed
in \refss{familycon} and let $\tilx_j\in \tile_j$ be the lift of $x_j\in e_j$. Then the trivialization of $\tilE$ defines isomorphisms
\[
    S_{z,j}:\,E_{x_j} \ \simeq \ \tilE_{\tilx_j} \ \to \ \CC^n, \qquad j=1\nek N, \ z\in \calO'.
\]
The isomorphisms $S_{z,j}$ depend on the trivialization of $\tilE$, i.e., on the connection $\n_z$. The direct sum $S_z= \bigoplus_jS_{z,j}$ is an
isomorphism $S_z:C^\b(K,E)\to C^\b(K,\gam(z))$ between the complex \refe{tilC(K,alp)} and \refe{C(K,alp)}. It follows from \refe{nz-nz0} that $S_z$
is complex differentiable at $z_0$, i.e., there exists a linear map $s:C^\b(K,E)\to \CC^{n\cdot N}$ such that
\eq{tilpa(z)}\notag
    S_z \ = \ S_{z_0}\ + \ (z-z_0)\,s\ + \ o(z-z_0).
\end{equation}

\subsection{The integration map}\label{SS:integr}
For each $z\in \calO'$ and for each $j=1\nek N$, the flat connection $\n_z$ defines an isomorphism $T_{j,z}:E|_{e_j}\to E_{x_j}\times e_j$. Thus, we
can define the {\em integration map}
\eq{integr}
    I_z:\,\Ome^\b(M,E)\ \to C^\b(K,E)
\end{equation}
by the formula
\eq{integr2}
    I_z(\ome) \ = \ \bigoplus_{1\le j\le N}\, \int_{e_j}\,T_{j,z}(\ome).
\end{equation}
By the de Rham theorem, $I_z$ is a morphism of complexes, i.e., $I_z\circ\n_z= \tilpa(z)\circ I_z$, which induces an isomorphism of cohomology. Also
it follows from \refe{nz-nz0} that $I_z$ is complex differentiable at $z_0$.

Finally, we consider the morphism of complexes
\eq{tilintegr}
    J_z \ := \ S_z\circ I_z\circ \psi(z):\, W^\b\ \to C^\b, \qquad z\in \calO'.
\end{equation}
This map is complex differentiable at $z_0$ and induces an isomorphism of cohomology.

\subsection{The cone complex}\label{SS:cone}
The cone complex $\Cone^\b(J_z)$ of the map $J_z$ is given by the sequence of vector spaces
\[
    \Cone^j(J_z) \ := \ W^j\oplus C^{j-1}\big(\,K,\gam(z)\,\big) \ \simeq \ \CC^{m_j}\oplus \CC^{n\cdot k_{j-1}}, \qquad j=0\nek d,
\]
with differentials
\[
    \hatpa_j(z) \ = \ \begin{pmatrix} d_j(z)&0\\J_{z,j}&\pa(\gam(z))\end{pmatrix},
\]
where $J_{z,j}$ denotes the restriction of $J_z$ to $W^j$. This is a family of acyclic complexes with differentials $\hatpa_j(z)$, which are complex
differentiable at $z_0$. The standard bases of $\CC^{m_j}\oplus \CC^{n\cdot k_{j-1}}$ define an element $c\in \Det(\Cone^\b(J_z))$ which is
independent of $z\in \calO'$. Using the isomorphism \refe{isomorphism}, we hence obtain for each $z\in \calO'$ the number $\phi_{\Cone^\b(J_z)}(c)\in
\CC\backslash\{0\}$. From the discussion in \refss{sketchnonacyclic} it follows that this number is equal to the ratio \refe{f(z)}. Hence, to finish
the proof of the \refp{holnearalp} it remains to show that the function $z\mapsto \phi_{\Cone^\b(J_z)}(c)$ is complex differentiable at $z_0$. This
follows immediately from \refl{asyclicanalytic}.

\section{Comparison between the Refined Analytic and the Farber-Turaev torsions}\label{S:comparison}

In this section we calculate the ratio of the refined analytic and the Farber-Turaev torsion. As a corollary, we conclude that the refined
analytic torsion is a holomorphic section on the whole space $\Rep$ and not only on the subset $\Rep\backslash\Sig(M)$ of smooth points.

First, we need to introduce some additional notations.

\subsection{The $\eta$-invariant}\Label{SS:etainv}
First, we recall the definition of the $\eta$-function of a non-self-adjoint elliptic operator $D$, cf. \cite{Gilkey84}. 
Let $D:C^\infty(M,E)\to C^\infty(M,E)$ be an elliptic differential operator of order $m\ge 1$ with self-adjoint leading symbol. Assume that $\tet$ is
an Agmon angle for $D$ (cf., for example, Definition~3.3 of \cite{BrKappelerRAT}). Let $\Pi_>$  (resp. $\Pi_<$) be a pseudo-differential projection
whose image contains the span of all generalized eigenvectors of $D$ corresponding to eigenvalues $\lam$ with $\RE\lam>0$ (resp. with $\RE\lam<0$)
and whose kernel contains the span of all generalized eigenvectors of $D$ corresponding to eigenvalues $\lam$ with $\RE\lam\le0$ (resp. with
$\RE\lam\ge0$). For all complex $s$ with $\RE{s}<-d/m$, we define the $\eta$-function of $D$ by the formula
\eq{eta}
    \eta_{\tet}(s,D) \ = \ \zet_\tet(s,\Pi_>,D) \ - \ \zet_\tet(s,\Pi_<,-D),
\end{equation}
where $\zet_\tet(s,\Pi_>,D):= \Tr(\Pi_>D^s)$ and, similarly, $\zet_\tet(s,\Pi_<,D):= \Tr(\Pi_<D^s)$. Note that, by definition, the purely imaginary
eigenvalues of $D$ do not contribute to $\eta_\tet(s,D)$.

It was shown by Gilkey, \cite{Gilkey84}, that $\eta_\tet(s,D)$ has a meromorphic extension to the whole complex plane $\CC$ with isolated simple
poles, and that it is regular at $0$. Moreover, the number $\eta_\tet(0,D)$ is independent of the Agmon angle $\tet$.

Since the leading symbol of $D$ is self-adjoint, the angles $\pm\pi/2$ are principal angles for $D$. Hence, there are at most finitely many
eigenvalues of $D$ on the imaginary axis. Let $m_+(D)$ (resp., $m_-(D)$) denote the number of eigenvalues of $D$, counted with their algebraic
multiplicities, on the positive (resp., negative) part of the imaginary axis. Let $m_0(D)$ denote the algebraic multiplicity of 0 as an eigenvalue of
$D$.

\defe{etainv}
The $\eta$-invariant $\eta(D)$ of\/ $D$ is defined by the formula
\eq{etainv}
    \eta(D) \ = \ \frac{\eta_\tet(0,D)+m_+(D)-m_-(D)+m_0(D)}2.
\end{equation}
\edefe
As $\eta_\tet(0,D)$ is independent of the choice of the Agmon angle $\tet$ for $D$, cf. \cite{Gilkey84}, so is $\eta(D)$.

\rem{etainv}
Note that our definition of $\eta(D)$ is slightly different from the one proposed by Gilkey in \cite{Gilkey84}. In fact, in our notation, Gilkey's
$\eta$-invariant is given by $\eta(D)+m_-(D)$. Hence, reduced modulo integers, the two definitions coincide. However, the number $e^{i\pi\eta(D)}$
will be multiplied by $(-1)^{m_-(D)}$ if we replace one definition by the other. In this sense, \refd{etainv} can be viewed as a {\em sign
refinement} of the definition given in \cite{Gilkey84}.
\erem

Let $\alp\in \Rep$ be a representation of the fundamental group of $M$ and let $E_\alp\to M$ be the vector bundle defined by $\alp$, cf.
\refss{Ealp}. We denote by $\na$ the flat connection on $\Ea$. Fix a Riemannian metric $g^M$ on $M$ and denote by
\eq{etaalp}
    \eta_\alp \ = \ \eta\big(\B_\even(\na,g^M)\big)
\end{equation}
the $\eta$-invariant of the corresponding odd signature operator $\B(\na,g^M)$, cf. \refd{oddsign}.

\subsection{The number $r_\ccomp$}\label{rccomp}
For every integer homology class $\xi\in H_1(M,\ZZ)$ and every $\alp\in$ \linebreak$\Rep$, we denote by $\det\alp(\xi)$ the determinant of the value
of $\alp$ on any closed curve $\gam$ representing $\xi$, $[\gam]=\xi$.

Let $L_{d-1}(p)\in H^{d-1}(M,\ZZ)$ denote the component in dimension $d-1$ of the Hirzebruch L-polynomial $L(p)$ in the Pontrjagin classes of $M$ and
let $\hatL_1\in H_1(M,\ZZ)$ denote the Poincar\'e dual of $L_{d-1}(p)$.

\lem{rccomp}
The function
\eq{rccomp0}
    \alp \ \mapsto \ r(\alp) \ := \ \big|\,\det\,\alp(\hatL_1)\,\big|^{1/2}\cdot e^{\pi\IM\eta_\alp} \ \in \ \RR_+
\end{equation}
is locally constant on $\Rep$. In particular, if $\ccomp\subset \Rep$ is a connected component of $\Rep$, which contains a unitary representation
$\alp_0$, then $\eta_{\alp_0}$ is real and $\big|\det\alp(\hatL_1)\big| = 1$, hence, $r(\alp)=1$ for all $\alp\in \ccomp$.
\elem
\prf
Following Farber \cite{Farber00AT} we denote by $\Arg_\alp$ the unique cohomology class in $H^1(M,\CC/\ZZ)$ such that for every closed curve
$\gam\in M$ we have
\eq{IArg}
    \det\big(\,\alp([\gam])\,\big) \ = \ \exp\big(\, 2\pi i\<\Arg_\alp,[\gam]\>\,\big),
\end{equation}
where $\<\cdot,\cdot\>$ denotes the natural pairing $H^1(M,\CC/\ZZ)\times{}H_1(M,\ZZ) \to  \CC/\ZZ$. Then
\[
    \log\, r(\alp) \ = \ \pi\, \IM\,\big(\, \eta_\alp \ - \ \<\Arg_\alp,\hatL_1\>\,\big).
\]
Suppose $\alp_t$ ($t\in[0,1]$) is a smooth family of representations. From Theorem~12.3 and Lemma~12.6 of \cite{BrKappelerRAT} we conclude that
\[
    \frac{d}{dt}\,\eta_{\alp_t} \ = \frac{d}{dt}\,\<\Arg_{\alp_t},\hatL_1\>.
\]
Hence, $\frac{d}{dt}r(\alp_t)= 0$.
\eprf
\defe{rccomp}
For each connected component $\ccomp\subset \Rep$ we denote by $r_\ccomp$ the value of the function $r$ on $\ccomp$.
\edefe
\refl{rccomp} implies that $r_\ccomp=1$ if $\ccomp$ contains a unitary representation.

\subsection{The homology class $\bet_\eps$}\label{SS:beteps}
We need the following
\lem{L-w}
Let $M$ be a closed oriented manifold of odd dimension $d=2n-1$. Let $L_{d-1}(p)\in H^{d-1}(M,\ZZ)$ denote the component in dimension $d-1$ of the
Hirzebruch L-polynomial $L(p)$ in the Pontrjagin classes of $M$. Then the reduction of $L_{d-1}(p)$ modulo 2 is equal to the $(d-1)$-Stiefel-Whitney
class $w_{d-1}(M)\in H^{d-1}(M,\ZZ_2)$  of\/ $M$.
\elem
\prf
For any homology class $\xi\in H_{d-1}(M,\ZZ)$ there exists a smooth oriented submanifold $X_\xi\subset M$, representing $\xi$. Then
$\<L_{d-1}(p),\xi\>$ is equal to the signature $\sig(X_\xi)$ of $X_\xi$. The parity of $\sig(X_\xi)$ is equal to the parity of the Euler
characteristic $\chi(X_\xi)$ of $X_\xi$, which, in turn, is equal to $\<w_{d-1}(M),X_\xi\> = \<w_{d-1}(X_\xi),X_\xi\>$. Thus we conclude that
\[
    \<L_{d-1}(p)-w_{d-1}(M),\xi\> = 0 \qquad \MOD\ 2,
\]
for any homology class $\xi\in H_{d-1}(M,\ZZ)$.
\eprf

We denote by $\hatL_1\in H_1(M,\ZZ)$ the Poincar\'e dual of $L_{d-1}(p)$ and by $c(\eps)\in H_1(M,\ZZ)$ the characteristic class of the Euler
structure $\eps$, cf. \cite{Turaev90} or Section~5.2 of \cite{FarberTuraev00}.
\cor{L-c}
The class $\hatL_1(p)+c(\eps)\in H_1(M,\ZZ)$ is divisible by 2, i.e. there exists a (not necessarily unique) homology class $\bet_\eps\in H_1(M,\ZZ)$
such that
\eq{L-c}
    -\,2\,\bet_\eps \ = \ \hatL_1(p)+c(\eps).
\end{equation}
\ecor
\prf
It is shown on page~209 of \cite{FarberTuraev00} that the reduction of $c(\eps)$ modulo 2 is equal to the Poincar\'e dual of the Stiefel-Whitney class $w_{d-1}(M)$.
Hence, it follows from \refl{L-w} that the reduction of $\hatL_1(p)+c(\eps)$ is the zero element of $H_1(M,\ZZ_2)$.
\eprf

The equality \refe{L-c} defines $\bet_\eps$ modulo two-torsion elements in $H_1(M,\ZZ)$. We fix a solution of \refe{L-c} and for the rest of the
paper $\bet_\eps$ denotes this solution.

\subsection{Comparison between the Farber-Turaev and the refined analytic torsions}\label{SS:Ray-Singer}
One of the main results of this paper is the following extension of the Cheeger-M\"uller theorem about the equality between the Reidemeister and the
Ray-Singer torsions:

\th{Ray-Singer}
Suppose $M$ is a closed oriented odd dimensional manifold. Let $\eps$ be an Euler structure on $M$ and let $\gro$ be a cohomological orientation of
$M$. Then, for each connected component $\ccomp$ of \/ $\Rep$, there exists a constant $\tet^\ccomp= \tet^\ccomp_\gro\in \RR/2\pi\ZZ$, depending on
$\gro$ (but not on $\eps$), such that,
\eq{tet-gro}
    \tet^\ccomp_{-\gro} \ \equiv \ \tet^\ccomp_\gro \ + \ n\pi, \qquad \MOD 2\pi,
\end{equation}
and for any representation $\alp\in \ccomp$,
\eq{logTTTur}
    \frac{\rat(\alp)}{\rtt(\alp)} \ = \ e^{i\tet^\ccomp_\gro}\cdot r_\ccomp\cdot   \det\alp(\bet_\eps),
\end{equation}
where $\bet_\eps\in H_1(M,\ZZ)$ is the homology class defined in \refe{L-c} and $r_\ccomp>0$ is defined in \refd{rccomp}. If the connected component
$\ccomp$ contains a unitary representation $\alp_0$, then $r_\ccomp=1$.
\eth
As an immediate corollary of \reft{Ray-Singer} we obtain
\cor{Ray-Singer}
If the representations $\alp_1, \alp_2$ belong to the same connected component of \linebreak$\Rep$ then
\eq{aa=tt}
    \frac{\rat(\alp_1)}{\rat(\alp_2)} \ = \ \frac{\rtt(\alp_1)}{\rtt(\alp_2)}\cdot \frac{\det\alp_1(\bet_\eps)}{\det\alp_2(\bet_\eps)}.
\end{equation}
\ecor
\subsection{Proof of \reft{holomsection}}\label{SS:prholomsection}
Before proving \reft{Ray-Singer} let us note that, since the right hand side of the equality \refe{logTTTur} is obviously holomorphic in $\alp$, it
follows from this equality that $\rat(\alp)$ is a holomorphic section of $\Dbundle$, cf. \refd{holsection}. Hence, \reft{holomsection} is proven.
\hfill$\square$

\subsection{Proof of \reft{Ray-Singer}}\label{SS:sketchRS}
In \refss{RSnormof} below, we use the calculations of the Ray-Singer norm of the Farber-Turaev torsion from \cite{FarberTuraev00} and the
calculation of the Ray-Singer norm of the refined analytic torsion from \cite{BrKappelerRATdetline} to compute the absolute value of the left
hand side of \refe{logTTTur}. More precisely we conclude that (cf. \refe{normofration2})
\eq{logTTTurabsv}
    \left|\, \det\alp(\bet_\eps)^{-1}\cdot\frac{\rat(\alp)}{\rtt(\alp)}\,\right| \ = \  r_\ccomp.
\end{equation}
By \refp{holnearalp}, ${\rat(\alp)}/{\rtt(\alp)}$ is an analytic function on the set $\Rep\backslash\Sig(M)$ of non-singular points of $\Rep$.
Further, $\det\alp(\bet_\eps)$ is obviously a polynomial function on $\Rep$.  Hence, the function
\[
    \alp \ \mapsto \ \det\alp(\bet_\eps)^{-1}\cdot\frac{\rat(\alp)}{\rtt(\alp)}
\]
is holomorphic on $\Rep\backslash\Sig(M)$. By \refe{logTTTurabsv} the absolute value of this function is locally constant $\Rep\backslash\Sig(M)$. It
follows that the function itself is locally constant, i.e., there exists a locally constant real valued function
$\tet_{\eps,\gro}:\Rep\backslash\Sig(M)\to \RR/2\pi\ZZ$ such that
\eq{logTTTur2}
    \frac{\rat(\alp)}{\rtt(\alp)} \ = \ e^{i\tet_{\eps,\gro}(\alp)}\cdot r_\ccomp\cdot   \det\alp(\bet_\eps), \qquad \alp\in \Rep\backslash\Sig(M).
\end{equation}

In \refl{continuity}, we show  that the function ${\rat(\alp)}/{\rtt(\alp)}$ is continuous on $\Rep$. Hence, $\tet_{\eps,\gro}(\alp)$ extends to a
continuous function on $\Rep$. Since $\tet_{\eps,\gro}$ is locally constant on the open dense subset $\Rep\backslash\Sig(M)$, which has only finitely
many connected components, it is also locally constant on $\Rep$. In other words, $\tet_{\eps,\gro}(\alp)$ depends only on the connected component
$\ccomp$ containing $\alp$.

To finish the proof of \reft{Ray-Singer} it remains to prove that $\tet_{\eps,\gro}$ is independent of $\eps$ and satisfies \refe{tet-gro}. This is
done in \refss{deponepsgro}.

\subsection{The Ray-Singer norm of the Farber-Turaev and the refined analytic torsions}\label{SS:RSnormof}
Let $\RSn{\cdot}$ denote the Ray-Singer norm on the determinant line $\Det(H^\b(M,E_\alp))$, cf.
\cite{Quillen85,BisZh92,FarberTuraev00,BrKappelerRATdetline}. Theorem~10.2 of \cite{FarberTuraev00} states that
\eq{RSnormrtt}
    \RSn{\rtt(\alp)} \ = \ \big|\,\det\alp(c(\eps))\,\big|^{1/2}.
\end{equation}
Further, by Theorem~11.3 of our previous paper \cite{BrKappelerRATdetline},
\eq{RSnormRAT}
    \RSn{\rat} \ = \ e^{\pi\,\IM \eta_\alp}.
\end{equation}
Combining \refe{RSnormrtt} and \refe{RSnormRAT} we obtain
\eq{normofration}
    \left|\,\frac{\rat(\alp)}{\rtt(\alp)}\,\right| \ = \ \big|\,\det\alp\big(c(\eps)\big)\,\big|^{-1/2}\cdot e^{\pi\,\IM \eta_\alp}.
\end{equation}

Since for any two homology classes $a,b\in H_1(M,\ZZ)$ we have $\det\alp(a+b)= \det\alp(a)\cdot\det\alp(b)$, in view of \refe{L-c} we obtain
\eq{alp(c)}
    \big|\, \det\alp\big(c(\eps)\big)\,\big| \ = \  \big|\, \det\alp\big(\,c(\eps)+\hatL_1\,)\,\big|\cdot \big|\, \det\alp(\hatL_1)\,\big|^{-1} \ = \
    \big|\, \det\alp(\bet_\eps)\,\big|^{-2}\cdot \big|\, \det\alp(\hatL_1)\,\big|^{-1}.
\end{equation}
Substituting \refe{alp(c)} into \refe{normofration} we obtain, using \refe{rccomp0},
\eq{normofration2}
 \begin{aligned}
    \left|\,\frac{\rat(\alp)}{\rtt(\alp)}\,\right| \ &= \ \big|\,\det\alp(\bet_\eps)\,\big|\cdot
            \Big(\, \big|\, \det\alp(\hatL_1)\,\big|^{1/2}\cdot e^{\pi\,\IM \eta_\alp}\,\Big)
    \\ &= \ \big|\,\det\alp(\bet_\eps)\,\big|\cdot r_\ccomp.
 \end{aligned}
\end{equation}

\lem{continuity}
The function \
 \(
    \alp \ \mapsto \ \frac{\rat(\alp)}{\rtt(\alp)}
  \) \
is continuous on $\Rep$.
\elem
\prf
The proof is similar (but easier) to the proof of \refp{holnearalp}. The only difference is that we now assume that $\alp_0\in \Rep$ is an arbitrary
(possibly singular) point and that the connection $\n_z$ depends merely continuously on $z$. Correspondingly, throughout the proof, one should
replace the words ``complex differentiable" by ``continuous".
\eprf

\subsection{Dependence of {}$\tet_{\eps,\gro}$ on the Euler structure and the cohomological orientation}\label{SS:deponepsgro}
From \refl{continuity} we conclude that $\tet_{\eps,\gro}$ is locally constant on $\Rep$. For each connected component $\ccomp$ of $\Rep$ denote
by $\tet_{\eps,\gro}^\ccomp$ the value of $\tet_{\eps,\gro}$ on $\ccomp$. To finish the proof of \reft{Ray-Singer} it remains to show that
$\tet_{\eps,\gro}^\ccomp$ is independent of $\eps$ and satisfies \refe{tet-gro}. In the case of acyclic representations the independence of
$\tet^\ccomp_{\eps,\gro}$ of $\eps$ was first established by R.-T.~Huang \cite{Huang06}.

Recall that the group $H_1(M,\ZZ)$ acts freely and transitively on the set $\operatorname{Eul}(M)$ of all Euler structures on $M$, cf.
\cite{Turaev90,FarberTuraev00}. Suppose $\eps_1,\eps_2\in \operatorname{Eul}(M)$ are two Euler structures and let $h\in H_1(M,\ZZ)$ be such that
\[
    \eps_2 \ = h+\eps_1,
\]
where $h+\eps_1$ denotes the action of $h$ on $\eps_1$. By formula (5.3) of \cite{FarberTuraev00}
\eq{cepa1eps2}
    c(\eps_2) \ = \ 2h \ + \ c(\eps_1) \ \in \ H_1(M,\ZZ).
\end{equation}
Further, by the first displayed formula on page~211 of \cite{FarberTuraev00}
\eq{rhoeps1eps2}
    \rho_{\eps_2,\gro}(\alp) \ = \ \det\alp(h)\cdot \rho_{\eps_1,\gro}(\alp).
\end{equation}
Combining \refe{cepa1eps2} and \refe{rhoeps1eps2} with \refe{L-c}, we conclude that
\eq{alpeps1eps2}
    \frac{\det\alp(\bet_{\eps_2})}{\det\alp(\bet_{\eps_1})} \ = \ \det\alp(\bet_{\eps_2}-\bet_{\eps_1}) \ = \ \det\alp(h)^{-1} \ = \
    \frac{\rho_{\eps_1,\gro}(\alp)}{\rho_{\eps_2,\gro}(\alp)}.
\end{equation}
Comparing \refe{alpeps1eps2} with \refe{logTTTur} we conclude that $\tet^\ccomp_{\eps,\gro}$ is independent of $\eps$.

It is shown in \S6.3 of \cite{FarberTuraev00} that
\eq{rhoeps-o}
    \rho_{\eps,-\gro}(\alp) \ = \ (-1)^n\, \rtt(\alp).
\end{equation}
Comparing this equality with \refe{logTTTur} we conclude that $e^{i\tet_{\eps,-\gro}^\ccomp}= (-1)^n\cdot e^{i\tet_{\eps,\gro}^\ccomp}$, which is
equivalent to \refe{tet-gro}. The proof of \reft{Ray-Singer} is now complete. \hfill$\square$

\subsection{Comparison with the Farber-Turaev absolute torsion}\label{SS:abstor}
An immediate application of \reft{Ray-Singer} concerns the notion of the {\em absolute torsion} introduced by Farber and Turaev in \cite{FarberTuraev99}. Suppose that
the Stiefel-Whitney class $w_{d-1}(M)\in H^{d-1}(M,\ZZ_2)$ vanishes, a condition always satisfied if $\dim{}M\equiv 3 (\MOD 4)$, cf. \cite{Massey60}. Then, by
\cite[\S3.2]{FarberTuraev99}, there exists an Euler structure $\eps$ such that $c(\eps)=0$. Assume, in addition, that the first Stiefel-Whitney class $w_1(E_\alp)$,
viewed as a homomorphism $H_1(M,\ZZ)\to \ZZ_2$, vanishes on the 2-torsion subgroup of $H_1(M,\ZZ)$. In this case there is also a canonical choice of the cohomological
orientation $\gro$, cf. \cite[\S3.3]{FarberTuraev99}. Then the Turaev torsion $\rtt(\alp)$ corresponding to any  $\eps$ with $c(\eps)=0$ and the canonically chosen
$\gro$ will be the same.

If the above assumptions on $w_{d-1}(M)$ and $w_1(E_\alp)$ are satisfied, then the number
\eq{abstor}
    \rabs(\alp) \ := \ \rtt(\alp)\ \in\ \CC, \qquad (\,c(\eps)\,=\,0\,),
\end{equation}
is canonically defined, i.e., is {\em independent of any choices}. It was introduced by Farber and Turaev, \cite{FarberTuraev99}, who called it the {\em absolute
torsion}.

In view of \refe{L-c} and the fact that $\hatL_1$ vanishes if $\dim{}M\equiv 3\ (\MOD 4)$, \reft{Ray-Singer} leads to the following
\cor{rat-abstor}
In addition to the assumptions made in \reft{Ray-Singer} suppose that $\dim{}M\equiv 3\ (\MOD 4)$ and that the 2-torsion subgroup of $H_1(M,\ZZ)$ is trivial. Then the
ratio $\rat(\alp)/\rabs(\alp)$ is a locally constant function on $\Rep$ and its absolute value is equal to 1. \footnote{Added in proof: Recently, Huang \cite{Huang06}
proved that if there exists a continuous path of representations, connecting $\alp$ with a unitary representation, then $\rat(\alp)/\rabs(\alp)= \pm
e^{-i\pi\rho_\alp}$, where $\rho_\alp=\eta_\alp-(\rank\alp)\eta_\trivial$ is the $\rho$-invariant of $E_\alp$.}
\ecor

\section{Application to the eta-invariant}\label{S:Fraberth}

As an application of \reft{Ray-Singer} we establish a relationship between the $\eta$-invariant and  the phase of the Farber-Turaev torsion
which improves and generalizes a theorem of Farber \cite{Farber00AT} and an earlier result of ours, cf. \refr{Farber} below.

\subsection{Phase of the Farber-Turaev torsion of a unitary representation}\Label{SS:argTTur}
Recall that if \linebreak$\alp\in \Repo$ is an acyclic representation, then we view the refined analytic torsion $\rat(\alp)$ as a non-zero complex
number, via the canonical isomorphism $\Det(H^\b(M,E_\alp))\simeq \CC$. We denote the phase of a complex number $z$ by $\Ph(z)\in [0,2\pi)$ so that
$z= |z|e^{i\Ph(z)}$.

\prop{extofFarber}
Suppose that $\alp_1,\, \alp_2\in \Repo$ are acyclic unitary representations which lie in the same connected component of\, $\Rep$. Then, modulo
$2\pi\,\ZZ$,
\meq{extofFarber2}
    \Ph\big(\,\rtt(\alp_1)\,\big)  + \pi\, \eta_{\alp_1}  +
    2\,\pi\,\big\<\Arg_{\alp_1},\bet_\eps\big\>
    \ \equiv \ \Ph\big(\,\rtt(\alp_2)\,\big)  + \pi\, \eta_{\alp_2} +
    2\,\pi\, \big\<\Arg_{\alp_2},\bet_\eps\big\>.
\end{multline}
\eprop
\prf
By formula (14.10) of \cite{BrKappelerRAT}, for any acyclic unitary representation $\alp$ we have
\eq{PhT0}
    \Ph\big(\,\rat(\alp)\,\big) \ = \ -\,\pi\,\eta_\alp \ + \ \pi\,(\rank \alp)\, \eta_\trivial,
    \qquad \MOD\, 2\pi\,\ZZ.
\end{equation}
Hence,
\eq{PhT}
    \Ph\big(\,\rat(\alp_1)\,\big) \ - \ \Ph\big(\,\rat(\alp_2)\,\big) \ = \ \pi\,\big(\, \eta_{\alp_2}-\eta_{\alp_1}\,\big)
    \qquad \MOD\, 2\pi\,\ZZ,
\end{equation}

From \refe{aa=tt} and \refe{IArg} we obtain, $\MOD 2\pi\ZZ$,
\eq{Ph-Ph}
 \begin{aligned}
    \Ph\big(\,\rat(\alp_1)\,\big) \ - \ \Ph\big(\,\rat(\alp_2)\,\big) \ &\equiv \
    \Ph\big(\,\rtt(\alp_1)\,\big) \ - \ \Ph\big(\,\rtt(\alp_2)\,\big)
    \\ &+ \ 2\,\pi\,\big\<\Arg_{\alp_1},\bet_\eps\big\> \ - \ 2\,\pi\,\big\<\Arg_{\alp_2},\bet_\eps\big\>.
 \end{aligned}
\end{equation}
Combining \refe{PhT} with \refe{Ph-Ph} we obtain \refe{extofFarber2}.
\eprf

\subsection{Sign of the absolute torsion}\Label{SS:signofAT}
Suppose that the Stiefel-Whitney class $w_{d-1}(M)=0$ and that the  first Stiefel-Whitney class $w_1(E_\alp)$, viewed as a homomorphism $H_1(M,\ZZ)\to \ZZ_2$,
vanishes on the 2-torsion subgroup of $H_1(M,\ZZ)$. Then the Farber-Turaev absolute torsion \refe{abstor} is defined. If $\alp\in \Repho$ is an acyclic unitary
representation, then $\rabs(\alp)$ is real, cf. Theorem~3.8 of \cite{FarberTuraev99} and, hence,
\[
    e^{i\Ph(\rabs(\alp))} \ = \ \sign\big(\,\rabs(\alp)\,\big).
\]
Note also, that since $c(\eps)=0$ it follows from \refe{L-c} that $2\bet_\eps= -\hatL_1$. Therefore,
\eq{halL1-bet}
    2\pi\, \<\Arg_\alp,\bet_\eps\> \ \equiv \ -\,\pi\, \<\Arg_\alp,\hatL_1\>, \qquad \MOD\ 2\pi\,\ZZ.
\end{equation}
Recall that  $\hatL_1$ vanishes if $\dim{M}\equiv 3$ ($\MOD\,4$).

From \refp{extofFarber} and \refe{halL1-bet} we now obtain the following 
\cor{Farber}
Suppose that $\alp_1,\, \alp_2\in \Repo$ are acyclic unitary representations which lie in the same connected component of\, $\Rep$. Suppose that he
first Stiefel-Whitney class $w_1(E_{\alp_1})= w_1(E_{\alp_2})$ vanishes on the 2-torsion subgroup of $H_1(M,\ZZ)$.
\begin{enumerate}
\item
If\/ $\dim{M}\equiv 3$ ($\MOD\,4$), then
\[
    \sign\big(\,\rabs(\alp_1)\,\big)\cdot e^{i\pi\eta_{\alp_1}} \ = \ \sign\big(\,\rabs(\alp_2)\,\big)\cdot e^{i\pi\eta_{\alp_2}}.
\]
\item
If\/ $\dim{M}\equiv 1$ ($\MOD\,4$) and $w_{d-1}(M)=0$, then
\[
    \sign\big(\,\rabs(\alp_1)\,\big)\cdot e^{i\pi( \eta_{\alp_1}- \< \Arg_{\alp_1},\hatL_1\>)}
    \ = \ \sign\big(\,\rabs(\alp_2)\,\big)\cdot e^{i\pi(\eta_{\alp_2}- \< \Arg_{\alp_2},\hatL_1\>)}.
\]
\end{enumerate}
\ecor

\rem{Farber}
For the special case when there is a real analytic path $\alp_t$ of {\em unitary} representations connecting $\alp_1$ and $\alp_2$ such that $\alp_t$ is acyclic for
all but
finitely many values of $t$, \refc{Farber} was established by Farber, using a completely different method,%
\footnote{Note that Farber's definition of the $\eta$-invariant differs from ours by a factor of 2. Moreover, the sign in front of
$\<\Arg_{\alp_1},\hatL_1\>$  in \cite{Farber00AT} should be the opposite one.} see \cite{Farber00AT}, Theorems~2.1 and 3.1. In
\cite[\S14.11]{BrKappelerRAT} we succeeded in eliminating the assumption of the existence of a real analytic path $\alp_t$ and assumed only that
the representations $\alp_1$ and $\alp_2$ lie in the same connected component of a certain subset of the set of acyclic representations.
\refc{Farber} improves on this result by showing that it is enough to assume that $\alp_1$ and $\alp_2$ lie in the same connected component of
$\Rep$.
\erem

\providecommand{\bysame}{\leavevmode\hbox to3em{\hrulefill}\thinspace} \providecommand{\MR}{\relax\ifhmode\unskip\space\fi MR }
\providecommand{\MRhref}[2]{%
  \href{http://www.ams.org/mathscinet-getitem?mr=#1}{#2}
} \providecommand{\href}[2]{#2}


\begin{thebibliography}{10}

\bibitem{APS1}
M.~F. Atiyah, V.~K. Patodi, and I.~M. Singer, \emph{Spectral asymmetry and
  {R}iemannian geometry. {I}}, Math. Proc. Cambridge Philos. Soc. \textbf{77}
  (1975), 43--69.

\bibitem{APS2}
\bysame, \emph{Spectral asymmetry and {R}iemannian geometry. {II}}, Math. Proc.
  Cambridge Philos. Soc. \textbf{78} (1975), no.~3, 405--432.

\bibitem{BeGeVe}
N.~Berline, E.~Getzler, and M.~Vergne, \emph{Heat kernels and {Dirac}
  operators}, Springer-Verlag, 1992.


\bibitem{BisZh92}
J.-M. Bismut and W.~Zhang, \emph{An extension of a theorem by {Cheeger} and
  {M\"uller}}, Ast\'erisque \textbf{205} (1992).

\bibitem{BrKappelerRAT}
M.~Braverman and T.~Kappeler, \emph{{Refined Analytic Torsion}},
  \texttt{arXiv:math.DG/0505537}.

\bibitem{BrKappelerRATdetline}
\bysame, \emph{{Refined Analytic Torsion as an Element of the Determinant
  Line}}, IHES preprint M/05/49, \texttt{arXiv:math.GT/0510532}.

\bibitem{BrKappelerRATshort}
\bysame, \emph{A refinement of the {Ray}-{Singer} torsion}, C.R. Acad. Sci.
  Paris \textbf{341} (2005), 497--502.

\bibitem{BrKappelerRATdetline_BH}
\bysame, \emph{{Comparison of the refined analytic and the
  Burghelea-Haller torsions}}, arXiv:math.DG/0606398.


\bibitem{BFK3}
D.~Burghelea, L.~Friedlander, and T.~Kappeler, \emph{Asymptotic expansion of
  the {Witten} deformation of the analytic torsion}, Journal of Funct. Anal.
  \textbf{137} (1996), 320--363.

\bibitem{BurgheleaHaller_Euler}
D.~Burghelea and S.~Haller, \emph{{{Euler} Structures, the Variety of
  Representations and the {Milnor}-{Turaev} Torsion}},
  \texttt{arXiv:math.DG/0310154}.

\bibitem{BurgheleaHaller_function}
\bysame, \emph{Torsion, as a function on the space of representations},
  \texttt{arXiv:math.DG/0507587}.

\bibitem{BurgheleaHaller_function2}
\bysame, \emph{{Complex valued Ray-Singer torsion}},
  arXiv:math.DG/0604484.

\bibitem{Farber00AT}
M.~Farber, \emph{Absolute torsion and eta-invariant}, Math. Z. \textbf{234}
  (2000), no.~2, 339--349.

\bibitem{FarberTuraev99}
M.~Farber and V.~Turaev, \emph{Absolute torsion}, Tel Aviv Topology Conference:
  Rothenberg Festschrift (1998), Contemp. Math., vol. 231, Amer. Math. Soc.,
  Providence, RI, 1999, pp.~73--85.

\bibitem{FarberTuraev00}
\bysame, \emph{Poincar\'e-{R}eidemeister metric, {E}uler structures, and
  torsion}, J. Reine Angew. Math. \textbf{520} (2000), 195--225.

\bibitem{Gilkey84}
P.~B. Gilkey, \emph{The eta invariant and secondary characteristic classes of
  locally flat bundles}, Algebraic and differential topology---global
  differential geometry, Teubner-Texte Math., vol.~70, Teubner, Leipzig, 1984,
  pp.~49--87.

\bibitem{GoldmanMillson88}
W.~Goldman and J.~Millson, \emph{The deformation theory of representations of
  fundamental groups of compact {K}\"ahler manifolds}, Inst. Hautes \'Etudes
  Sci. Publ. Math. (1988), no.~67, 43--96.

\bibitem{HormanderSCV}
L.~H{\"o}rmander, \emph{An introduction to complex analysis in several
  variables}, third ed., North-Holland Mathematical Library, vol.~7,
  North-Holland Publishing Co., Amsterdam, 1990.

\bibitem{Huang06}
R.-T. Huang, \emph{{Refined analytic torsion of lens spaces}},
  \texttt{arXiv:math.DG/0602231}.

\bibitem{Massey60}
W.~S. Massey, \emph{On the {S}tiefel-{W}hitney classes of a manifold}, Amer. J.
  Math. \textbf{82} (1960), 92--102.

\bibitem{Quillen85}
D.~Quillen, \emph{Determinants of {Cauchy}-{Riemann} operators over a {Riemann}
  surface}, Funct. Anal. Appl. \textbf{14} (1985), 31--34.

\bibitem{Rudyak_book}
Y.~B. Rudyak, \emph{On {T}hom spectra, orientability, and cobordism}, Springer
  Monographs in Mathematics, Springer-Verlag, Berlin, 1998, With a foreword by
  Haynes Miller.

\bibitem{Turaev01}
V.~G. Turaev, \emph{Introduction to combinatorial torsions}, Lectures in Mathematics
  ETH Z\"urich, Birkh\"auser Verlag, Basel, 2001, Notes taken by Felix Schlenk.

\bibitem{Turaev86}
\bysame, \emph{{R}eidemeister torsion in knot theory}, Russian Math.
  Survey \textbf{41} (1986), 119--182.

\bibitem{Turaev90}
\bysame, \emph{Euler structures, nonsingular vector fields, and
  {R}eidemeister-type torsions}, Math. USSR Izvestia \textbf{34} (1990),
  627--662.

\bibitem{Wall60}
C.~T.~C. Wall, \emph{Determination of the cobordism ring}, Ann. of Math. (2)
  \textbf{72} (1960), 292--311.

\end{thebibliography}
\end{document}